\documentclass[11pt]{amsart}
\usepackage{amsmath,amssymb, graphicx, amscd,latexsym,here}
\makeatletter
\newtheorem{Theorem}{Theorem}
\newtheorem{Lemma}[Theorem]{Lemma}
\newtheorem{Corollary}[Theorem]{Corollary}
\newtheorem{Proposition}[Theorem]{Proposition}

\newcommand\la{\lambda}
\newcommand\vphi{\varphi}

\newcommand\al{\alpha}
\newcommand\La{\Lambda}
\newcommand\si{\sigma}
\newcommand\be{\beta}

\newcommand\ga{\gamma}

\newcommand\cS{\mathcal  S}

\newcommand\BC{ {\mathbb C}}
\newcommand\BN{ {\mathbb  N}}
\newcommand\BQ{ {\mathbb  Q}}
\newcommand\BZ{{\mathbb  Z}}
\newcommand\BR{ {\mathbb  R}}

\newcommand\bfv{\mbox {\bf  v}}

\newcommand\bfx{\mbox {\bf  x}}

\newcommand\bfw{\mbox {\bf  w}}
\newcommand\bfm{\mbox {\bf  m}}
\newcommand\bfn{\mbox {\bf  n}}

\newcommand\bfz{\mbox {\bf  z}}
\newcommand\bfy{\mbox {\bf  y}}
\newcommand\bfa{\mbox {\bf  a}}
\newcommand\bfb{\mbox {\bf  b}}

\newcommand\nl{\newline}

\newcommand\GL{\rm{GL}\/}
\newcommand\mini{\rm{min}\/}

\newcommand\id{\rm{id}}

\newcommand\mv{\rm{m-vector}}

\newcommand\inv{^{-1}}

\def\mapright#1{\smash{\mathop{\longrightarrow}\limits^{{#1}}}}

\def\mapdown#1{\Big\downarrow\rlap{$\vcenter{\hbox{$#1$}}$}}

\def\inv{^{-1}}

\begin{document}
\title[ Topology of  polar weighted homogeneous hypersurfaces
]
{Topology of  polar weighted homogeneous hypersurfaces
}

\author
[M. Oka ]
{Mutsuo Oka }
\address{\vtop{
\hbox{Department of Mathematics}
\hbox{Tokyo  University of Science}
\hbox{26 Wakamiya-cho, Shinjuku-ku}
\hbox{Tokyo 162-8601}
\hbox{\rm{E-mail}: {\rm oka@rs.kagu.tus.ac.jp}}
}}
\keywords {Polar weighted homogeneous, Polar action}
\subjclass[2000]{14J17, 32S25}

\begin{abstract}
Polar weighted homogeneous polynomials are the class of special polynomials 
of real variable $x_i,y_i,\, i=1,\dots,n$ with $z_i=x_i+\sqrt{-1}y_i$
 which enjoys a `` polar action''. In many aspects, their
behavior looks like that of complex weighted homogeneous polynomials.
We study basic properties of hypersurfaces which are defined by
polar weighted homogeneous polynomials. 
\end{abstract}
\maketitle

\maketitle

\section{Introduction}
We consider a polynomial
$f(\bfz,\bar \bfz)=\sum_{\nu,\mu}c_{\nu,\mu}\bfz^\nu\bar\bfz^\mu$
where $\bfz=(z_1,\dots,z_n)$, $\bar \bfz=(\bar z_1,\dots,\bar z_n)$,
$ \bfz^\nu=z_1^{\nu_1}\cdots z_n^{\nu_n}$ for $ \nu=(\nu_1,\dots,\nu_n)$
(respectively
$\bar\bfz^\mu=\bar z_1^{\mu_1}\cdots\bar z_n^{\mu_n}$  for 
$ \mu=(\mu_1,\dots,\mu_n)$
as usual. Here $\bar z_i$ is the complex conjugate of $z_i$. 
Writing $z_i=x_i+\sqrt{-1}y_i$, it is easy to see that $f$ is a
polynomial of $2n$-variables $x_1,y_1,\dots, x_n,y_n$.
 Thus $f$ can be understood as an real analytic function
$f:\,\BC^n\to \BC$.
 We call $f$ {\em a mixed polynomial}
of
$z_1,\dots, z_n$.

A mixed polynomial
$f(\bfz,\bar \bfz)$ is called {\em polar weighted homogeneous}
if there exists   integers
$q_1,\dots, q_n$ and $p_1,\dots, p_n$ and non-zero integers $m_r,\,m_p$
such that 
\begin{eqnarray*}
&\gcd(q_1,\dots,q_n)=1,\quad\gcd(p_1,\dots, p_n)=1,\, \\
&\sum_{i=1}^n q_j( \nu_j+\mu_j)=m_r,\quad
\sum_{i=1}^n p_j (\nu_j-\mu_j)=m_p,\quad \text{if}\,\, c_{\nu,\mu}\ne 0
\end{eqnarray*}
We say $f(\bfz,\bar\bfz)$ is {\em a polar weighted homogeneous of 
radial weight type\nl $(q_1, \dots, q_n ; m_r) $ and of polar weight type
$(p_1,\dots, p_n;m_p)$}.
We define vectors of rational numbers
 $(u_1,\dots, u_n)$ and $(v_1,\dots, v_n)$ by $u_i=q_i/m_r,\mv_i=p_i/m_p$
and we call them {\em the normalized radial
 (respectively  polar) weights}.
Using a polar coordinate
$(r,\eta)$ of $\BC^*$ where $r>0$ and $\eta\in S^1$ with
$S^1=\{\eta\in \BC\,|\, |\eta|=1\}$, we define {\em a polar $\BC^*$-action} on
$\BC^n$ by
\begin{eqnarray*}
& (r,\eta)\circ \bfz=(r^{q_1}\eta^{p_1}z_1,\dots, 
r^{q_n}\eta^{p_n}z_n),\quad
(r,\eta)\in \BR^+\times S^1\\
&(r,\eta)\circ\bar \bfz=\overline{(r,\eta)\circ \bfz}=(r^{q_1}\eta^{-p_1}\bar z_1,\dots, 
r^{q_n}\eta^{-p_n}\bar z_n).
\end{eqnarray*}
Then 
$f$ satisfies the
functional equality
\begin{eqnarray}\label{polar-weight}
f((r,\eta)\circ (\bfz,\bar \bfz) )=r^{m_r}\eta^{m_p}f(\bfz,\bar\bfz).\end{eqnarray}
This notion is introduced by  Ruas-Seade-Verjovsky \cite{R-S-V}
implicitly
and then by Cisneros-Molina \cite{Molina}.

It is easy to see that such a polynomial defines a global fibration
\[
 f:\BC^n-f\inv(0)\to \BC^{*}.
\]
The purpose of this paper is to study the topology of the hypersurface
$F=f\inv(1)$ for a given polar weighted homogeneous polynomial, which is
a fiber of the above fibration. Note that 
F has a canonical stratification
\[
 F=\amalg_{I\subset \{1,2,\dots,n\}} F^{*I},\quad F^{*I}=F\cap \BC^{*I}
\]
Our main result is Theorem \ref{main-result}, which describes the
topology of $F^{*I}$ for a  simplicial polar weighted polynomial.

 \section{Polar weighted homogeneous hypersurface}
This section is  the preparation for the later sections.
Proposition 2 and Proposition 3 are added for consistency but they 
 are essentially
known in the series of works by J. Seade
and coauthors \cite{R-S-V,Seade,Pichon-Seade,Pichon-Seade2,SeadeBook}.
\subsection{Smoothness of a mixed hypersurface}
Let $f(\bfz,\bar\bfz)$ be a mixed polynomial and we consider a
hypersurface
$V=\{\bfz\in \BC^n;f(\bfz,\bar\bfz)=0\}$.
Put $z_j=x_j+ iy_j$. Then
$f(\bfz,\bar\bfz)$ is  a real analytic function of $2n$ variables
$(\bfx,\bfy)$ with $\bfx=(x_1,\dots,x_n)$ and $\bfy=(y_1,\dots, y_n)$.
Put $f(\bfz,\bar\bfz)=g(\bfx,\bfy)+i\,h(\bfx,\bfy)$ where $g,\,h$ are real
analytic functions.
Recall that
\begin{eqnarray*}
&\frac{\partial}{\partial z_j}=\frac 12 \left(
\frac{\partial}{\partial x_j}-i\frac{\partial}{\partial y_j}\right),
\quad\frac{\partial}{\partial \bar z_j}=\frac 12 \left(
\frac{\partial}{\partial x_j}+i\frac{\partial}{\partial y_j}
\right)\end{eqnarray*}
Thus 
\begin{eqnarray*}
&\frac{\partial k}{\partial z_j}=\frac 12
\left(\frac{\partial k}{\partial x_j}-i\frac{\partial k}{\partial
 y_j}\right),\quad \frac{\partial k}{\partial \bar z_j}=\frac 12
\left(\frac{\partial h}{\partial x_j}+i\frac{\partial k}{\partial
 y_j}\right)
\end{eqnarray*}
for any analytic function $k(\bfx,\bfy)$.
Thus for a complex valued function $f$, we define
\begin{eqnarray*}
&\frac{\partial f}{\partial z_j}=\frac{\partial g}{\partial z_j}+i\frac{\partial h}{\partial
 z_j},
\quad
\frac{\partial f}{\partial \bar z_j}=
\frac{\partial g}{\partial \bar z_j}+i\frac{\partial g}{\partial
 \bar z_j}
\end{eqnarray*}

We assume that $g,\,h$ are non-constant polynomials. Then $V$ is real
codimension two subvariety. Put 
\begin{eqnarray*}
& d_{\BR}g(\bfx,\bfy)=(\frac{\partial g}{\partial x_1},\dots,\frac{\partial
 g}{\partial x_n},\frac{\partial g}{\partial y_1},\dots, 
\frac{\partial g}{\partial y_n})\in \BR^{2n}
\\
& d_{\BR}h(\bfx,\bfy)=(\frac{\partial h}{\partial x_1},\dots,\frac{\partial
 h}{\partial x_n},\frac{\partial h}{\partial y_1},\dots, 
\frac{\partial h}{\partial y_n})\in \BR^{2n}
\end{eqnarray*}
For a complex valued mixed polynomial, we use the notation:
\begin{eqnarray*}
&d f(\bfz,\bar\bfz)=(\frac{\partial f}{\partial z_1},\dots,\frac{\partial
 f}{\partial z_n})\in \BC^n,\quad
\bar d f(\bfz,\bar\bfz)=(\frac{\partial f }{\partial\bar z_1},\dots,\frac{\partial
 f}{\partial\bar z_n})\in \BC^n
\end{eqnarray*}
Recall that a point $\bfz\in V$ is a singular point of $V$ if and only
if
two vectors
$dg(\bfx,\bfy),\,dh(\bfx,\bfy)$ are linearly dependent over
$\BR$ (see Milnor \cite{Milnor}). This condition is not so easy to be checked, as the calculation
of $g(\bfx,\bfy),\,h(\bfx,\bfy)$ from a given $f(\bfz,\bar\bfz)$ is not
immediate. However we have
\begin{Proposition}\label{Singular-condition}
The following two conditions are equivalent.
\begin{enumerate}
\item  $\bfz\in V$ is a singular point of $V$ and $\dim_{\BR}(V,\bfz)=2n-2$.
\item There exists a complex 
number $\alpha,\,|\al|=1$ such that 
$ \overline{df(\bfz,\bar\bfz)}=\al\, \bar d f(\bfz,\bar\bfz)$.
\end{enumerate}
\end{Proposition}
\begin{proof} First assume that $d_{\BR}g,\,d_{\BR}h$ are linearly
 dependent at $\bfz$.
Suppose for example that $dg(\bfx,\bfy)\ne 0$ and write
$dh(\bfx,\bfy)=t\,dg(\bfx,\bfy)$ for some $t\in \BR$.
This implies that 
\begin{eqnarray*}
& \frac {\partial f}{\partial x_j}=(1+ti)\frac{\partial g}{\partial x_j},
\quad
 \frac {\partial f}{\partial y_j}=(1+ti)\frac{\partial g}{\partial y_j},
\quad\text{thus}\,\,\\
&\frac {\partial f}{\partial z_j}=(1+ti)\left(\frac{\partial g}{\partial x_j}
-i\frac{\partial g}{\partial y_j}\right),\,\,
\frac {\partial f}{\partial \bar z_j}=(1+ti)\left(\frac{\partial g}{\partial x_j}
+i\frac{\partial g}{\partial y_j}\right).
\end{eqnarray*}
Thus 
\begin{eqnarray*}
&d  f(\bfz,\bar\bfz)=(1+ti)
\left(\frac{\partial g}{\partial x_1}-i\frac{\partial g}{\partial y_1},
\dots, \frac{\partial g}{\partial x_n}-i\frac{\partial g}{\partial y_n}\right)=
2(1+ti) d_{\bfz} g(\bfz,\bar\bfz)
\\
&\bar d f(\bfz,\bar\bfz)=(1+ti)
\left(\frac{\partial g}{\partial x_1}+i\frac{\partial g}{\partial y_1},
\dots, \frac{\partial g}{\partial x_n}+i\frac{\partial g}{\partial y_n}\right)=
2(1+ti) d_{\bar\bfz} g(\bfz,\bar\bfz)
\end{eqnarray*}
Here $d_{\bfz}g=(\frac{\partial g}{\partial z_1},\dots,\frac{\partial
 g}{\partial z_n})$ and 
$d_{\bar\bfz}g=(\frac{\partial g}{\partial \bar z_1},\dots,\frac{\partial
 g}{\partial \bar z_n})$.
As $g$ is a real valued polynomial, using the equality $\overline{d_{\bfz}g(\bfx,\bfy)}=d_{\bar\bfz}g(\bfx,\bfy)$ we get 
\[
 \overline{df(\bfz,\bar\bfz)}=\frac{1-ti}{1+ti}\bar df(\bfz,\bar\bfz).
\]
Thus it is enough to take $\al=\frac{1-ti}{1+ti}$.

Conversely assume that $\overline{df(\bfz,\bar\bfz)}=\al \bar
 df\bfz,\bar\bfz$
for some $\al=a+bi$ with $a^2+b^2=1$.
Using the notations
\[
 d_x g=(\frac{\partial g}{\partial x_1},\dots,\frac{\partial g}{\partial
 x_n}),
\quad
d_y g=(\frac{\partial g}{\partial y_1},\dots,\frac{\partial g}{\partial
 y_n}),\, \text{etc},
\]
we get\begin{eqnarray*}
&(1-a)d_xg +b \,d_y g=-b\, d_xh-(1+a)d_y h\\
&-b\, d_xg\,+\,(1-a)d_y g=(a+1)d_x h-b\, d_y h.
\end{eqnarray*}
Solving these equations, we get
\[
d_{\BR}\,g=(d_xg,d_yg)=\frac{-2b}{(1-a)^2+b^2}d_{\BR}\, h
\]
which proves the assertion.
\end{proof}
\subsection{Polar weighted homogeneous hypersurfaces}
Let $f$ be a polar weighted homogeneous polynomial
of radial weight type $(q_1,\dots, q_n;m_r)$
and of polar weight type $(p_1,\dots, p_n;m_p)$. By differentiating 
(\ref{polar-weight}) in \S 1, we get
\begin{eqnarray}
&m_r f(\bfz,\bar\bfz)=\sum_{i=1}^n q_i(\frac{\partial f}{\partial
 z_i}z_i+\frac{\partial f}{\partial
 \bar z_i}\bar z_i
)\\
&m_p f(\bfz,\bar\bfz)=\sum_{i=1}^n p_i(\frac{\partial f}{\partial
 z_i}z_i-\frac{\partial f}{\partial
 \bar z_i}\bar z_i
).
\end{eqnarray}
We call these equalities {\em Euler equalities.}
Recall that $\BC^n$ has the canonical hermitian inner product defined by
\[
 (\bfz,\bfw)=z_1\bar w_1+\cdots+z_n\bar w_n.
\]
Identifying $\BC^n$ with $\BR^{2n}$ by $\bfz\longleftrightarrow (\bfx,\bfy)$,
the Euclidean inner product of $\BR^{2n}$ is given as
$(\bfz,\bfw)_{\BR}=\Re (\bfz,\bfw)$. Or we can also write  as
\[
 (\bfz,\bfw)_{\BR}=\frac 12\left( (\bfz,\bfw)+(\bar\bfz,\bar\bfw)\right).
\]
\begin{Proposition}
For any $\alpha\ne 0$, the fiber $F_\al:=f\inv(\al)$ 
is a smooth $2(n-1)$ real-dimensional manifold and it is canonically
 diffeomorphic to $F_1=f\inv(1)$.
\end{Proposition}
\begin{proof}
Take a point $\bfz\in F_\al$. We consider two particular vectors
$\bfv_r,\,\bfv_\theta\in T_{\bfz} \BC^n$ which are the tangent vectors of the
respective orbits
of $\BR$ and $S^1$:
\begin{eqnarray*}
&\bfv_r=\frac{d(r\circ \bfz)}{dr}|_{r=1}=(q_1z_1,\dots, q_nz_n),\quad \\
&\bfv_\theta=\frac{d(e^{i\theta}\circ \bfz)}{d\theta}|_{\theta=0}=
(ip_1 z_1,\dots, ip_n z_n).
\end{eqnarray*}
Taking the differential of the equality
\[
 f((r,\exp(i\theta))\circ \bfz))=r^{m_r}\exp(m_p\theta
 i)f(\bfz,\bar\bfz),
\]
we see that
$df_z:T_{\bfz}\BC^n\to T_\al \BC^*$ 
satisfies
\[
 df_z(\bfv_r)=m_r |\al|\frac{\partial}{\partial r},
\quad
df_z(\bfv_\theta)=m_p \frac{\partial}{\partial \theta}
\]
where $(r,\theta)$ is the polar coordinate of $\BC^*$.
This implies that $f:\BC^n\to \BC$ is a submersion at $\bfz$. Thus 
$F_\al$ is a smooth codimension 2 submanifold.
A diffeomorphism $\vphi_\al: F_1\to F_\al$ is simply given as
$ \vphi(\bfz)=(r^{1/m_r},\exp^{i\theta/m_p})\circ \bfz$
where $\al=r\exp(i\theta)$.
\end{proof}
The above proof does not work for $\al=0$. 
Recall that the polar $\BR^+$ action along the radial direction is  written in real coordinates as
$$r\circ(\bfx,\bfy)=(
r^{q_1}x_1,\dots,r^{q_n}x_n,r^{q_1}y_1,\dots,r^{q_n}y_n ),\quad
r\in \BR^+.$$

\begin{Proposition}\label{smoothness}
Let $V=f\inv(0)$. Assume that $q_j>0$ for any $j$.
 Then $V$ is contractible to the origin $O$.
If further $O$ is an isolated singularity of $V$,
$V \backslash\{O\}$ is smooth.
\end{Proposition}
\begin{proof}
 A canonical deformation retract $\be_t:V\to V$ is given as
 $\be_t(\bfz)=t\circ \bfz$, $0\le t\le 1$.
(More precisely $\be_0(\bfz)=\lim_{t\to 0}\be_t(\bfz)$.) Then $\be_1=\id_V$ and $\be_0$ is the
 contraction to $O$.
Assume that $\bfz\in V\backslash\{O\}$ is a singular point.
Consider the decomposition into real analytic functions
  $f(z)=g(\bfx,\bfy)+ih(\bfx,\bfy)$.
Using the radial $\BR^+$-action, we see that
\begin{eqnarray}\label{Euler-real}
 g(r\circ (\bfx,\bfy))=r^{m_r}g(\bfx,\bfy),\quad
h(r\circ (\bfx,\bfy))=r^{m_r}h(\bfx,\bfy).
\end{eqnarray}
This implies that $g(\bfx,\bfy),\,h(\bfx,\bfy)$ are weighted homogeneous
 polynomials
of $(\bfx,\bfy)$ and the Euler equality can be restated as 
\begin{eqnarray*}
 m_r \, g(\bfx,\bfy) &=\sum_{j=1}^n p_j\left(x_j\frac{\partial g}{\partial x_j}(\bfx,\bfy)+
y_j \frac{\partial g}{\partial y_j}(\bfx,\bfy)\right)\\
m_r\, h(\bfx,\bfy)&=\sum_{j=1}^n p_j\left(x_j\frac{\partial h}{\partial x_j}(\bfx,\bfy)+
y_j \frac{\partial h}{\partial y_j}(\bfx,\bfy)\right).
\end{eqnarray*}
Differentiating the equalities (\ref{Euler-real}) in $r$, we get
\[
 \frac{\partial g}{\partial x_j}(r\circ
 (\bfx,\bfy))=r^{m_r-q_j}\frac{\partial g}{\partial x_j}(\bfx,\bfy),\quad
\frac{\partial h}{\partial x_j}(r\circ (\bfx,\bfy))
 =r^{m_r-q_j}\frac{\partial h}{\partial x_j}(\bfx,\bfy).
\]
 This  implies that these differentials are also  weighted
 homogeneous polynomials of degree $m_r-q_j$.
Thus the jacobian matrix
\[
 \left(\frac{\partial (g,h)}{\partial{(x_i,y_i)}}(r\circ
 (\bfx,\bfy))
\right)
\]
is the same with the jacobian matrix at $\bfz=(\bfx,\bfy)$ up to scalar
multiplications in the column vectors by
$r^{m_r-q_1},\dots, r^{m_r-q_n},r^{m_r-q_1},\dots, r^{m_r-q_n}$ respectively.
Thus any points of the orbit
$r\circ (\bfx,\bfy),\,r>0$ are singular points of $V$.
This is a contradiction to the assumption that $O$ is an isolated
 singular point of $V$, as $\lim_{r\to 0}r\circ (\bfx,\bfy)=O$.
\end{proof}
\begin{Proposition} (Transversality)\label{Transversality}
Under the same assumption as in Proposition \ref{smoothness}, the
 sphere
 $S_\tau=\{\bfz\in\BC^n; |\bfz|=\tau\}$
intersects transversely with $V$ for any $\tau>0$.
\end{Proposition}
\begin{proof}
Let $\phi(\bfx,\bfy)=\|\bfz\|^2=\sum_{j=1}^n(x_j^2+y_j^2)$.
Then $S_\tau$ intersects transversely with $V$ if and only if
the gradient vectors
$d_{\BR}g,d_{\BR}h,d_{\BR}\phi$ are linearly  independent over $\BR$.
Note that $d_{\BR}\phi(\bfx,\bfy)=2(\bfx,\bfy)$.
Suppose that 
 the sphere $S_{\|\bfz\|}$ is tangent to $V$ at $\bfz=(\bfx,\bfy)\in V$.
Then we have for example, a linear relation
$dg(\bfx,\bfy)=\al\, dh(\bfx,\bfy)+\be\, d\phi(\bfx,\bfy)$ with some $\al,\be\in
 \BR$.
Note that the tangent vector $\bfv_r$ to the $\BR^+$-oribit
 is tangent to $V$ and it is written
$\bfv_r=(q_1x_1,\dots,q_n x_n,q_1 y_1,\dots, q_n y_n)$ as a real vector.
Then we have
\begin{eqnarray*}
0=&\frac{dg(r\circ (\bfx,\bfy)}{dr}|_{r=1}=
\sum_{j=1}^n q_j\left(x_j\frac{\partial g}{\partial x_j}(\bfx,\bfy)+
y_j \frac{\partial g}{\partial y_j}(\bfx,\bfy)\right)\\
&=(\bfv_r (\bfx,\bfy),dg(\bfx,\bfy))_{\BR}
 \\
&=(\bfv_r (\bfx,\bfy),\al\, dh (\bfx,\bfy))+(\bfv_r (\bfx,\bfy),\be\, d\phi
 (\bfx,\bfy))_{\BR}\\
&=2\be\sum_{j=1}^nq_j(x_j^2+y_j^2)
\end{eqnarray*}
as $(\bfv_r (\bfx,\bfy),dh(\bfx,\bfy))_{\BR}=0$ by the same reason.
This is  the case only if $\be=0$ which is impossible as $V\backslash \{O\}$
is non-singular by Proposition \ref{smoothness}.
\end{proof}

\subsubsection{ Remark} Let $f(\bfz,\bar\bfz)$ be a polar weighted homogeneous
 polynomial with respective weights
 $(q_1,\dots,q_n;m_r)$ and $(p_1,\dots,p_n;m_p)$.
Proposition \ref{smoothness} does not hold if the radial
weights
contain some negative $q_j$. Assume that $q_j\ge 0$ for any $j$ and
$I_0:=\{j|q_j=0\}$ is not empty. Then it is easy to
see that $f$ does not have monomial which does not contain any $z_i$
with $i\notin I_0$, as if such monomial exists, its radial degree is 0.
This implies that $V=f\inv(0)$ contains  the coordinate subspace
$\BC^{I_0}=\{\bfz|z_i=0,\,i\notin I_0\}$.
We call $\BC^{I_0} $  {\em the canonical retract coordinate subspace}.
Then Proposition \ref{smoothness} can be modified as
{\em $\BC^{I_0}$ is a deformation retract of $V$.}
Of course, $\BC^{I_0}$ can be contracted to $O$ but this contraction is
not through the action and not related to the geometry of $V$.
\subsubsection{Example}
Consider the following examples.
 \begin{eqnarray*}
&g_1(\bfz,\bar \bfz)
 =z_1^{a_1}\bar z_2+\cdots+z_n^{a_n}\bar z_1,\,\forall a_i\ge
 1,\,\text{and}\,\, a_j\ge 2\, (\exists j)\\
&g_2(\bfz,\bar \bfz)=z_1^{a_1}\bar z_2+\cdots+z_{n-1}^{a_{n-1}}\bar z_n
+z_n^{a_n},\,\forall a_i\ge 1.
\end{eqnarray*}
\begin{Proposition}
\begin{enumerate}
\item 
 The radial
      weight vector
$(q_1,\dots, q_n)$ of $g_1(\bfz,\bar\bfz)$ is 
     semi- positive, i.e. $q_j\ge 0$
for any $j$ if $a_i\ge 1$ for any $i$.
($\exists j,\,a_j\ge 2$ is necessary for the existence of polar action.)
It is not strictly positive  if and only if $n=2m$ is even and
  either
(a) $a_1=a_3=\cdots=a_{2m-1}=1$ or (b) $a_2=a_4=\cdots=a_{2m}=0$.

In case (a) (respectively (b)), we have $q_2=q_4=\cdots=q_{2m}=0$
and $q_{2j+1}\ge 1,\,0\le j\le m-1$ (resp. $q_1=q_3=\cdots=q_{2m-1}=0$
and $q_{2j}\ge 1,\,1\le j\le m$ ).
\item 
 The radial    weight vector 
$(q_1,\dots, q_n)$ of $g_2(\bfz,\bar\bfz)$ is 
     semi-positive.
It is not strictly positive  if and only if
  $a_n=1$. Let $s$ be the integer such that
 $a_n=a_{n-2}=\cdots=a_{n-2s}=1$ and
      $a_{n-2s-2}\ge 2$.
Then $q_{n-1}=\dots=q_{n-2s+1}=0$ and $q_j\ge 1$ otherwise.
\end{enumerate}
\end{Proposition}
\begin{proof}
We first consider $g_1(\bfz,\bar\bfz)=z_1^{a_1}\bar z_2+\cdots+z_n^{a_n}\bar z_1$.
By an easy calculation, using the notation
 $a_{i+n}=a_i$ the normalized radial weigts $(u_1,\dots,u_n)$ are given as
\begin{eqnarray*}
& u_{j}=\frac{1}{a_1\cdots
 a_n-1}\sum_{i=0}^{m-1}(a_{j+2i+1}-1)a_{j+2i+2}\cdots a_{j+n-1},
\,n=2m\\
& u_{j}=\frac{1}{a_1\cdots a_n+1}\left(1+
\sum_{i=0}^{m-1}(a_{j+2i+1}-1)a_{j+2i+2}\cdots a_{j+n-1}\right),
\, n=2m+1\\
\end{eqnarray*}
and the assertion follows immediately from this expression.

Next we consider $g_2(\bfz,\bar\bfz)=z_1^{a_1}\bar
 z_2+\cdots+z_n^{a_{n-1}}\bar z_n+z_n^{a_n}$. Then 
the normalized radial weigts $(u_1,\dots,u_n)$ are given as
\begin{eqnarray*}
&u_j=\frac 1{a_j}-\frac 1{a_{j}a_{j+1}}+\cdots+(-1)^{n-j}\frac
 1{a_ja_{j+1}\cdots a_n}\\
&=\begin{cases}
&\frac{a_{j+1}-1}{a_{j}a_{j+1}}+\cdots+\frac {a_n-1}{a_ja_{j+1}\cdots a_n}
,\,
n-j:\, \text{odd}\\
&\frac{a_{j+1}-1}{a_{j}a_{j+1}}+\cdots+\frac {a_{n-1}-1}{a_ja_{j+1}\cdots
  a_{n-1}}
+\frac 1{a_ja_{j+1}\cdots a_n}\\
& \qquad\qquad\qquad n-j: \text{even}
\end{cases}
\end{eqnarray*}
As $a_i\ge 1$, the assertion follows from the above expression.
\end{proof}

\subsection{Simplicial mixed polynomial}
Let $f(\bfz,\bar\bfz)=\sum_{j=1}^{s}c_j\,\bfz^{{\bfn_j}}\bar\bfz^{\bfm_j}$
 be a mixed polynomial. Here we assume that $c_1,\dots, c_s\ne 0$.
Put 
\[
 \hat f(\bfw):=\sum_{j=1}^s c_j\,\bfw^{\bfn_j-\bfm_j}.
\]
We call $\hat f$ the {\em the associated Laurent polynomial.}
This polynomial plays an important role for the determination of the 
topology of the hypersurface $F=f\inv(1)$.
Note that
\begin{Proposition} If $f(\bfz,\bar\bfz)$ is a polar weighted 
homogeneous polynomial of polar weight type  $(p_1,\dots,p_n;m_p)$,
$\hat f(\bfw)$ is also a weighted homogeneous Laurent polynomial
of type $(p_1,\dots,p_n;m_p)$ in the complex variables $w_1,\dots,w_n$.
\end{Proposition}

A mixed polynomial $f(\bfz,\bar\bfz)$
is called {\em simplicial} if
the  exponent vectors $\{{\bfn_j}\pm \bfm_j\,|\,j=1,\dots, a\}$ 
are linearly independent in $\BZ^n$ respectively.
In particular, simplicity implies that $s\le n$. When $s=n$, we say that 
$f$ is {\em full.}
Put
$\bfn_j=(n_{j,1},\dots,n_{j,n})$,
$\bfm_j=(m_{j,1},\dots,m_{j,n})$ in $\BN^n$.
Assume that $s\le n$. 
Consider two integral matrix
$N=(n_{i,j})$ and
$M=(m_{i,j})$
where the $k$-th row vectors are $\bfn_k,\,\bfm_k$
respectively.
\begin{Lemma} Let $f(\bfz,\bar\bfz)$ be a  mixed polynomial as above.
If $f(\bfz,\bar\bfz)$ is simplicial, then $f(\bfz,\bar\bfz)$ is a polar weighted homogeneous polynomial.
In the case $s=n$, $f(\bfz,\bar\bfz)$ is simplicial if and only if
 $\det(N \pm M)\ne 0$.
\end{Lemma}
\begin{proof}
First we assume that $s=n$ and consider the system of linear equations
\begin{eqnarray}
&\label{radial}\begin{cases}
&(n_{1,1}+m_{1,1})u_1+\cdots+ (n_{1,n}+m_{1,n} )u_n=1\\
&\qquad\cdots\\
&(n_{n,1}+m_{n,1}) u_1+\cdots+ (n_{n,n}+m_{n,n} )u_n=1\\
\end{cases}
\\
&\label{polar}\begin{cases}
&(n_{1,1}-m_{1,1}) v_1+\cdots+ (n_{1,n}-m_{1,n}) v_n=1\\
&\qquad \cdots\\
&( n_{n,1}-m_{n,1}) v_1+\cdots+ (n_{n,n}-m_{n,n}   ) v_n=1\\
\end{cases}
\end{eqnarray}
It is easy to see that equations $(\ref{radial})$ and  $(\ref{polar})$
have solutions if
$\det\,N\pm M\ne 0$ which is equivalent for $f$ to be simplicial by definition.
Note that the solutions $(u_1,\dots, u_n)$ and $(v_1,\dots, v_n)$
are rational numbers. We call them {\em the normalized radial
 (respectively  polar) weights}.
Now let $m_r,\,m_p$ be the least common multiple of
the denominators of $u_1,\dots, u_n$ and $v_1,\dots, v_n$ respectively.
Then the weights are given as 
$q_j=u_jm_r ,\,p_j=v_j m_p,\, j=1,\dots, n$ respectively.

Now suppose that $s<n$. It is easy to choose positive integral vectors
${\bfn_j},\,j=s+1,\dots, n$ (and put $\bfm_j=0,\,j=s+1,\dots,n$)
such that 
$\det(\tilde N\pm \tilde M)\ne 0$,
where
$\tilde N$ and $\tilde M$ are $n\times n$-matrices adding $(n-s)$
row vectors $\bfn_{s+1},\dots,\bfn_n$. Then the assertion
follows from the case $s=n$.
This corresponds to considering the mixed polynomial:
\[
 f(\bfz,\bar\bfz)=\sum_{j=1}^s c_j \bfz^{\bfn_j}\bar\bfz_j^{\bfm_j}
+ 0\times \sum_{j=s+1}^n \bfz^{\bfn_j}.
\]
\end{proof}

\subsubsection{Example}
Let 
\begin{eqnarray*}
& f_{\bfa,\bfb}(\bfz,\bar \bfz)=z_1^{a_1}\bar z_2^{b_1}
+\dots+z_n^{a_n}\bar z_1^{b_n},\,a_i,b_i\ge 1,\,\forall i\\
&k(\bfz,\bar\bfz)=z_1^d(\bar z_1+\bar z_2)+\cdots+z_n^{d}(\bar z_n+\bar
 z_1),\,d\ge 2.
\end{eqnarray*}
The associated Laurent polynomials are
\begin{eqnarray*}
&\widehat{ f_{\bfa,\bfb}}(\bfw)=w_1^{a_1} w_2^{-b_1}
+\dots+w_n^{a_n} w_1^{-b_n}\\
&\hat k(\bfw)=w_1^d(1/ w_1+1/ w_2)+\cdots+w_n^{d}(1/ w_n+1/ w_1).
\end{eqnarray*}
\begin{Corollary}\label{polar-condition}
For  the polynomial  $f_{\bfa,\bfb}$, the following conditions are
 equivalent.
\begin{enumerate}
\item  $f_{\bfa,\bfb}$ is simplicial. 
\item $f_{\bfa,\bfb}$ is a polar weighted homogeneous polynomial.
\item ${\rm (SC)}\quad a_1\cdots a_n\ne b_1\cdots b_n.$
\end{enumerate}
\end{Corollary}
{\em Proof.}  The assertion follows from the equality:
\begin{eqnarray*}
&\det(\bfn\pm\bfm)=
 \det\,
\left(
\begin{matrix}
a_1&0&\cdots&\pm b_n\\
\pm b_1&a_2&\cdots&0\\
\vdots&\ddots&\ddots&\vdots\\
0&\cdots&\pm b_{n-1}&a_n
\end{matrix}
\right)\\
&
=\begin{cases}
&a_1a_2\dots a_n+(-1)^{n-1}b_1b_2\dots b_n \,\,\text{for}\, \bfn+\bfm\\
&\,a_1a_2\dots a_n-b_1b_2\dots b_n\,\,\,\text{for}\,\,\bfn-\bfm.
\end{cases}
\end{eqnarray*}
The polynomial $k(\bfz,\bar\bfz)$ is a polar weighted homogeneous 
polynomial with respective weight types
$(1,\dots,1;d+1)$ and $(1,\dots,1;d-1)$. However it is not simplicial.\qed

 Now we consider an example which does not satisfy the simplicial
 condition (SC) of
 Corollary
\ref{polar-condition}: 
$\phi_a:=z_1^a\bar z_1^a+\dots+z_n^a\bar z_n^a$.
This does not have any polar action
as they are polynomials of $|z_1|^2,\dots, |z_n|^2$ and it takes only
 non-negative values.
Note also that $\phi_a\inv(1)$ is real codimension 1 as
 $\phi_a(\bfx,\bfy)=\sum_{j=1}^n (x_j^2+y^2)^a$.

As a typical simplicial polar weighted polynomial, we consider
again the following two polar weighted polynomials.
 \begin{eqnarray*}
&g_1(\bfz,\bar \bfz)
 =z_1^{a_1}\bar z_2+\cdots+z_n^{a_n}\bar z_1,\,\forall a_i\ge
 1,\,\text{and}\,\, a_j\ge 2\, (\exists j)\\
&g_2(\bfz,\bar \bfz)=z_1^{a_1}\bar z_2+\cdots+z_{n-1}^{a_{n-1}}\bar z_n
+z_n^{a_n},\,\forall a_i\ge 1.
\end{eqnarray*}
The polynomial $g_1(\bfz,\bar \bfz)$ with $a_i\ge 2,\,(\forall i)$
is  a special case of $\sigma$-twisted  Brieskorn
 polynomial  and 
has been studied intensively (\cite{R-S-V}). In out case, we only assume $a_i\ge 2$
for some $i$. 
 The existence of $i$ with
$a_i\ge 2$ is the condition for
the existence of polar action.
 We consider two hypersurfaces defined by
$V_i=g_i\inv(0)$ for $i=1,2$.
The  condition for a hypersurface defined by a polar weighted
 homogeneous polynomial to have an isolated singularity
is more complicated than that of the
singularity defined by a complex anaytic hypersurface.
 For the above examples, we assert the following.
\begin{Proposition} \label{isolatedness}
For $V_1,\,V_2$,  we have the following criterion.
\begin{enumerate}
\item $V_i\cap \BC^{*n},\,i=1,2$ are non-singular.
\item $V_1=g_1\inv(0)$ has no singularity outside of the origin if and only if one of
      the following conditions is satisfied.
\begin{enumerate}
\item $n$ is odd.
\item $n$ is even and there are (at least) two indices $i,j\,(i<j)$
such that $a_i,\,a_j\ge 2$ and $j-i$ is odd.
\end{enumerate}
\item $V_2=g_2\inv(0)$
 has  no  singularity outside of the origin
if and only if 
 one of
      the following conditions is satisfied.
\begin{enumerate}
\item $a_n\ge 2$.
\item $a_n=1$, $n=2m+1$ is odd and $a_{2j-1}=1$ for any $1\le j\le m+1$.
\end{enumerate}
\end{enumerate}
\end{Proposition}
\begin{proof} We use Proposition \ref{Singular-condition}.
So assume that 
\[
( \sharp):\quad \overline{df(\bfz,\bar\bfz)}=\al\bar df(\bfz,\bar\bfz),
\quad |\al|=1.
\]

\vspace{.2cm}\noindent
(1) We consider $V_1$. Suppose $\bfz\in V_1\cap \BC^{*n}$ is a singular
 point.
Note that 
\[df_{\bfz}=(a_1z_1^{a_1-1}\bar z_2,\cdots, a_{n}z_{n}^{a_{n}-1}\bar
 z_1),\quad
 \bar df(\bfz,\bar\bfz)=(z_n^{a_n},z_1^{a_1},\dots, z_{n-1}^{a_{n-1}})
\]
$(\sharp)$ implies that
\begin{eqnarray}\label{eq1}
  a_j\bar z_j^{a_j-1}\bar z_{j+1}=\al z_{j-1}^{a_{j-1}},\,\,j=1,\dots,n,
\,\exists \al\in S^1. 
\end{eqnarray}
In this case, indices should be understood to be integers modulo $n$. So $z_{n+1}=z_1$,
 and so on.
 If $\bfz\in \BC^{*n}$, multiplying the absolute values of the both
 side of the above equality,
 this give a contradiction:
 $\prod_{i=1}^n a_i  |z_i|^{a_i}=\prod_{i=1}^n |z_i|^{a_i}$.

Now we consider the smoothness on $V_1\backslash\{O\}$.
Assume that $\bfz$ is a singular point of $V_1\backslash\{O\}$.
For simplicity, we may assume that $a_n\ge 2$ as $g_1$ is symmetric with
 the permutation $i\to i+1$. 

Assume that $z_\iota\ne 0$. 
Then the $(\iota+1)$-th component of $\bar df$ is non-zero.
Thus by $(\sharp)$, $(\iota+1)$-th component of $ df$ is also non-zero.
That is,  $z_{\iota+1}^{a_\iota-1}\bar z_{\iota+2}\ne 0$.
In particular, $z_{\iota+2}\ne 0$.  We repeat the same argument and
get a sequence of non-zero components
$z_\iota,z_{\iota+2},\dots$. Thus we arrive to the
 conclusion that either $z_{n-1}\ne 0$ (if $n-\iota$ is odd ) or
$z_{n}\ne 0$ ( if $n-\iota$ is even ).

--If $n-\iota$ is odd and $z_{n-1}\ne 0$,
 the last component of $df$ is non-zero and we
 have $z_n,z_1\ne 0$ as we have assumed that $a_n\ge 2$. 
This creates two non-zero sequence $z_n,z_2,z_4,\dots$ and $z_1,z_3,\dots$.
Thus 
we conclude that $\bfz\in \BC^{*n}$,
which is impossible by the first argument.

--If $n-\iota$ is even, $z_\iota,z_{\iota+2},\dots,z_{n}\ne 0$.
Thus we see that the first component of $\bar df$
is non-zero. By the same argument, we get 
non-zero sequence $z_2,z_4,\dots$.

Thus to show that $\bfz\in \BC^{*n}$, it is enough to show that
 $z_{n-1}\ne 0$.

\noindent
(a) Assume first $n$ is odd.
If $\iota$ is even, then we see that
$z_\iota,z_{\iota+2},\dots, z_{n-1}\ne 0$ and we are done.

If $\iota$ is odd, we get $z_n\ne 0$, which implies the first component
of $\bar df(\bfz,\bar\bfz)$ is non-zero.
Thus as the second round, we have non-zero sequence $z_2,z_4,\dots$
which contains $z_{n-1}$. Thus we are done.

\noindent
(b)
Now we assume that $n$ is even but there is another integer $1\le i<n$
such that $a_i\ge 2$ and $a_n\ge 2$ and $i$ is odd.
If $\iota$ is odd,  we have shown that 
$\bfz\in \BC^{*n}$.

If $\iota$ is even, we get $z_n\ne 0$ and thus $z_2\ne 0$.
Then the sequence $z_2,z_4,\dots$ contains $z_{i-1}$. As $a_i\ge 2$,
looking at the $i$-th component of $df$, we get
$z_{i}\cdot z_{i+1}\ne 0$. Thus  we get a non-zero sequence
$z_{i}, z_{i+2},\dots$ which contains $z_{n-1}$, and we are done.

Now to show that one of  the conditions (a) or (b) is 
necessary,
we  assume that $n$ is even  and $a_{\nu}=1$ for any odd $\nu$ and 
$a_n\ge 2$. Thus putting $n=2m$,
 \[
  f=(z_1\bar z_2+z_2^{a_2}\bar z_3)+\cdots+(z_{2m-1}\bar z_{2m}+z_{2m}^{a_{2m}}\bar z_1).
\]
Consider the subvariety $z_1=z_3=\cdots=z_{n-1}=0$.
Then 
\[
 df(\bfz,\bar\bfz)=(\bar z_2,0,\bar z_4,0,\dots, \bar z_{2m},0),
\quad
\bar df(\bfz,\bar\bfz)=(z_n^{a_n},0,\dots, z_{2m-2}^{a_{2m-2}},0)
\]the condition $(\sharp)$ is written as
\[
 (\sharp)\quad
 z_2=\al\,z_n^{a_n},\, z_4=\al\,z_2^{a_2},\cdots,  z_{2m}=\al\,z_{2m-2}^{a_{2m-2}}
\]
which has real one-dimensional solution
\begin{eqnarray*}
&z_{2j}=\al^{\be_j}u^{\ga_j}\,(\quad j=1,\dots, m),\,
\quad \al^{\be_m}\, u^{\ga_m a_{2m}-1}=1\\
&\be_j=1+\sum_{i=1}^{j-1}a_{2(j-1)}a_{2(j-2)}\cdots a_{2(j-i)},\,
\ga_j=a_{2}a_4\cdots a_{2(j-1)}
\end{eqnarray*}

\noindent
(2) We consider the case $V_2$.
We will see first   $V_2\cap \BC^{*n}$ is non-singular.
Take a singular point of $V_2$. Then  we have some $\al\in S^1$ so that
 \[
 ( \sharp):\,\overline{df}(\bfz,\bar\bfz)=\al
 {\bar df}(\bfz,\bar\bfz).
\]
As we have
\begin{eqnarray*}
 &df(\bfz,\bar\bfz)=(a_1z_1^{a_1-1}\bar z_2,\cdots,
 a_{n-1}z_{n-1}^{a_{n-1}-1}\bar z_n,a_nz_n^{a_n-1}),\\
&\bar df(\bfz,\bar\bfz)=(0,z_1^{a_1},\dots, z_{n-1}^{a_{n-1}})
\end{eqnarray*}
we see that $(\sharp)$ implies that $z_1^{a_1-1}\bar z_2=0$. Thus there
 are no singularities on $V_2\cap \BC^{*n}$.
Suppose that $z_\iota\ne 0$ for some $\iota$.
If $\iota<n-1$, this implies $(\iota+1)$-th component of $\bar df(\bfz,\bar \bfz)$
is non-zero. Thus $(\sharp)$ implies that $(\iota+1)$-th component
of $df$ is non-zero.
 In particular,
$z_{\iota+2}$ is non-zero.
(Of course, $z_{\iota+1}\ne 0$ if $a_{\iota+1}>1$.) Repeating this argument, 
we arrive to the conclusion: either $z_{n-1}$ or $z_n$ is non zero.

First assume that $a_n\ge 2$. Comparing the last components of $df(\bfz,\bar \bfz)$
and $\bar d f(\bfz,\bar \bfz)$, we  observe that 
$z_{n-1}$ and $z_n$ are both non-zero.
Now we go in  the reverse direction.
As the $(n-1)$-th component 
of  $df$ is non-zero,
the corresponding $(n-1)$-th component $z_{n-2}^{a_{n-2}}$ of $\bar df(\bfz,\bar
 \bfz)$
is non-zero. Then $(n-2)$-th component of $df(\bfz,\bar \bfz)$
is non-zero.
Going downword, we see that $z\in \BC^{*n}$. However this is impossible,
 as we have already  seen above.

Next we assume that $a_n=1$ and $n$ is odd and 
 $a_{2j-1}=1$ for any $j$.
 Note that the last component of
 $df(\bfz,\bar\bfz)$ is 1. Thus $z_{n-1}\ne 0$. If $z_{n}\ne 0$, we
get a contradiction as above $\bfz\in \BC^{*n}$. Thus we may assume that $z_n=0$.
Comparing $(2j)$-components of $df(\bfz,\bar\bfz)$ and $\al\bar df(\bfz,\bar\bfz)$,
 we get
\[
 z_2=0,\, z_4=\al z_2^{a_2},\dots, z_{n-1}=z_{n-3}^{a_{n-3}}
\]
 which has no solution with $z_{n-1}\ne 0$.

 Now we show that the condition $(a)$ or $(b)$ in (3) is necessary.

 (i) Assume that $a_n=1$ and $n$ is even ans put $n=2m$. 
Let $s$ be the maximal integer
 such that $a_{2s}\ge 2$. If there does not exists
 such $s$, we put $s=0$.
 Non-isolated singularities are given  by the solutions of 
 \begin{eqnarray*}
  &z_2=z_4=\cdots=z_{2m}=0,\, z_{2j-1}=0,\,j \le s\\
  &z_{2s+3}=\al z_{2s+1}^{a_{2s+1}},\dots, z_{2m-1}=\al
   z_{2m-3}^{a_{2m-3}},\,
   1=\al z_{2m-1}^{a_{2m-1}}.
   \end{eqnarray*}

 (ii) Assume that $a_n=1$, $n=2m+1$  is odd, and there exists odd index such
 that  $a_{2j+1}\ge 2$.
 Put $s$ be the maximum integer of such $j$.
  Non-isolated singularities are given  by the solutions of
 \begin{eqnarray*}
&  z_1=z_3=\cdots=z_{2m+1}=0,\, z_{2j}=0,\,j\le s\\
  & z_{2s+4}=\al z_{2s+2}^{a_{2s+2}},\dots,
   z_{2m}=\al z_{2m-2}^{a_{2m}},\, 1=\al z_{2m}^{a_{2m}}.
   \end{eqnarray*}
\end{proof}
\subsubsection{Remark}
1. The polynomial $g_1(\bfz,\bar\bfz)=z_1^{a_1}\bar
 z_2+\cdots+z_n^{a_n}\bar z_1$ is an example of so-called {\em 
$\si$-twisted
 Brieskorn polynomial} if $a_i\ge 2,\, i=1,\dots,n$.
Let $\si$ be a permutation of $\{1,2,\dots,n\}$.
Then $\si$-twisted  Brieskorn polynomial is defined as
\[
 f_\si(\bfz,\bar\bfz)=z_1^{a_1}\bar z_{\si(1)}+\cdots+z_n^{a_n}\bar
 z_{\si(n)},\quad
a_1,\dots, a_n\ge 2.
\] and the
 corresponding
 assertions in Proposition 3 and 4 are proved in \cite{Seade}.
See also \cite{SeadeBook} for more systematical treatment for 
real analytic polynomials which define  Milnor fibrations.
In \cite{H-Lopez}, similar conditions for the isolatedness condition as
Proposition \ref{isolatedness} are considered.
For our purpose,  we call $f_\si(\bfz,\bar\bfz)$ a {\em weak
 $\si$-twisted
Brieskorn polynomial} if $\si\in \cS_n$ and $a_i\ge 1,\,\forall i$.

\vspace{.3cm}\noindent
2. Consider a product
$\BC^n=\BC^s\times \BC^{n-s}$ and  use  variables $\bfv\in\BC^s$
and $\bfw\in \BC^{n-s}$. Assume that there exist mixed polynomials
$h(\bfv,\bar\bfv)$ and $k(\bfw,\bar\bfw)$ so that
$f(\bfz,\bar\bfz)=h(\bfv,\bar\bfv)+k(\bfw,\bar\bfw)$.
$f(\bfz,\bar\bfz)$ is a polar weighted polynomial
if and only if $h(\bfv,\bar\bfv),\,k(\bfw,\bar\bfw)$ are polar weighted
polynomial
and it is known that $f\inv(1)$ is homotopic to the join $h\inv(1)\star
k\inv(1)$ if $f$ is polar weighted.
 Such a polynomial is called a
{\em polynomial of join type} (\cite{Molina}, see also \cite{OkaJoin}).

Now consider a weak $\si$-twisted Brieskorn polynomial
 $f_\si(\bfz,\bar\bfz)$.
If $\si$ has order $n$, it is (up to a change of ordering)
equal to the cyclic permutation $\si=(1,2,\dots,n)$
and $f_\si=g_1$. In general, $\si$ can be written as a product of
 mutually commuting cyclic permutations $\si=\tau_1\tau_2\cdots \tau_\nu$.
 Put $|\tau_i|=\{j|\tau_i(j)\ne j\}$ and put $f_{\tau_i}$ be the partial
 sum of monomials in $f(\bfz,\bar\bfz)$ written in variables $\{z_j|j\in
 |\tau_i|\}$.
 Thus $f_\si$ is a join type polynomial of $\nu$ weak $\tau_i$-twisted
 Brieskorn polynomial $f_{\tau_i}$.  Thus
 $f_\si(\bfz,\bar\bfz)$ has 
 an isolated singularity if and only if each polynomial
$f_{\tau_i}$ has an isolated singularity.
 A similar assertion is also proved in \cite{H-Lopez}.

\vspace{.3cm}
\noindent
3. Observe that the singularities of $V_1,\,V_2$ are on the canonical 
retract coordinate subspaces $\BC^{I_0}$. Note also that the polar
action is trivial on $\BC^{I_0}$. 

\subsection{Milnor fibration}
Let $f(\bfz,\bar\bfz)$ be a polar weighted homogeneous polynomial of 
radial weight  type $(q_1,\dots, q_n;m_r)$ and of polar weight
  type
\nl $(p_1,\dots,p_n;m_p)$.
Then 
\[
 f:\BC^n-f\inv(0)\to \BC^*
\] is a locally trivial fibration.
The local triviality is given by the action. 
In particular, the monodromy map
$h:F\to F$ is given by
$h(\bfz)=\exp(2\pi i/m_p)\circ \bfz=(z_1 \exp(2p_1\pi i/m_p),\dots,
z_n\exp(2p_n\pi i/m_p))$ where $F=f\inv(1)$ (\cite{R-S-V,Molina}).
\section{Topology of simplicial polar weighted homogeneous hypersurface}
Let $f(\bfz,\bar\bfz)=\sum_{j=1}^{s}c_j\,\bfz^{{\bfn_j}}\bar\bfz^{{\bfm_j}}$
 be a   polar weighted homogeneous polynomial of 
radial weight type $(q_1,\dots, q_n;m_r)$ and of polar  weight
type \nl
$(p_1,\dots,p_n;m_p)$. Let $F=f\inv(1)$ be the fiber.
\subsection{Canonical stratification of $F$ and the topology of each  stratum}
For any subset $I\subset \{1,2,\dots,n\}$, we define
\begin{eqnarray*}
&\BC^I=\{\bfz\,|\, z_j= 0,\,\, j\notin I\},
\,\, \BC^{*I}=\{\bfz\,|\, z_i\ne 0\,\,\text{iff}\, \,i\in I\},\,
\BC^{*n}=\BC^{*\{1,\dots,n\}}
\end{eqnarray*}
 and  we define mixed polynomials  $f^I$ by the restriction:
$f^I=f|_{\BC^{I}}$. For simplicity, we write
a point of $\BC^I$ as $\bfz_I$.
Put
$F^{*I}=\BC^{*I}\cap F$.
Note that $F^{*I} $ is a  non-empty subset of $\BC^{*I}$
if and only if $f^I(\bfz_I,\bar \bfz_I)$ is not constantly zero.
Now we observe that the hypersurface $F=f\inv(1)$
has the canonical stratification
\[
 F=\amalg_{I}\,F^{*I}.
\]
Thus it is essential to determine the topology of each stratum $F^{*I}$.
 Put $F^*:=F\cap \BC^{*n}$,
the open dense stratum and put $\hat F^*:=\hat f\inv(1)\cap \BC^{*n}$
where $\hat f(\bfw)$
is the associated Laurent weighted homogeneous 
polynomial.
\begin{Theorem}\label{main-result}
Assume  that $f(\bfz,\bar\bfz)$ is a simplicial polar weighted
 homogeneous polynomial and let $\hat f(\bfw)$
be the associated Laurent weighted homogeneous 
polynomial. Then there exists a canonical diffeomorphism
$\vphi:\BC^{*n}\to \BC^{*n}$ which 
gives an isomorphism of tiwo Milnor fibrations defined by
$f(\bfz,\bar\bfz)$ and $\hat f(\bfw)$:
\[
 \begin{matrix}
\BC^{*n}-f\inv(0)&\mapright{f}& \BC^*\\
\mapdown{\vphi}&& \mapdown{\id}\\
\BC^{*n}-{\hat f}\inv(0)&\mapright{\hat f}&\BC^*
\end{matrix}
\]
and it satisfies
$\vphi(F^{*n})=\hat
 F^{*n}$
and $\vphi$ is compatible with the respective canonical monodromy maps.
\end{Theorem}
\begin{proof}
Assume first that $s=n$ for simplicity.
Recall that
\[
 \hat f(\bfw)=\sum_{j=1}^n \,c_j \bfw^{\bfn_j-\bfm_j}.
\]
Let $\bfw=(w_1,\dots, w_n)$ be the complex coordinates of $\BC^n$ which
 is the ambient space of $\hat F$.
We construct $\vphi:\BC^{*n}\to \BC^{*n}$ so that $\vphi(\bfz)=\bfw$ satisfies
$$\bfw(\vphi(\bfz))^{\bfn_j-\bfm_j}=\bfz^{\bfn_j}\bar\bfz^{\bfm_j},
\quad\text{thus}\quad
\hat f(\vphi(\bfz))=f(\bfz).
$$
For the construction of $\vphi$, we use the polar coordinates
$(\rho_j,\theta_j)$ for $z_j\in \BC^*$ and polar coordinates
$(\xi_j,\eta_j)$ for $\bfw_j$.
Thus $\bfz_j=\rho_j\,\exp(i\theta_j)$ and $\bfw_j=\xi_j\,\exp(i\eta_j)$.
First we take $\eta_j=\theta_j$.
 Put
$\bfn_j=(n_{j,1},\dots,n_{j,n})$,
$\bfm_j=(m_{j,1},\dots,m_{j,n})$ in $\BN^n$.
Consider two integral matrix
$N=(n_{i,j})$ and
$M=(m_{i,j})$
where the $k$-th row vector are $\bfn_k,\,\bfm_k$
respectively.
Now taking the logarithm of
 the equality $\bfz^{\bfn_j}\bar\bfz^{\bfm_j}=\bfw^{\bfn_j-\bfm_j}$,
we get an  equivalent equality:
\begin{eqnarray*}
&(n_{j1}+m_{j1})\log \rho_1+\cdots+(n_{jn}+m_{jn})\log\rho_n\qquad\qquad\\
\quad &=
(n_{j1}-m_{j1})\log \xi_1+\cdots+(n_{jn}-m_{jn})\log\xi_n
\end{eqnarray*}
for $j=1,\dots, n$. This can be written as
\begin{eqnarray}
(N+M)\left(
\begin{matrix}
\log\rho_1\\\vdots\\\log\rho_n
\end{matrix}\right)=
(N-M)\left(
\begin{matrix}
\log\xi_1\\\vdots\\\log\xi_n
\end{matrix}\right)
\end{eqnarray}
Put  $(N-M)\inv (N+M)=(\la_{ij})\in \GL(n,\BQ)$.
Now we define $\vphi$ as follows.
\begin{multline*}
\vphi:\BC^{*n}\to\BC^{*n},\,\,\bfz=(\rho_1\exp(i\theta_1),\dots,\rho_n\exp(i\theta_n))\mapsto\\
\bfw=(\xi_1\exp(i\theta_1),\dots,\xi_n\exp(i\theta_n))
\end{multline*}
where $\xi_j$ is given by
$ \xi_j=\exp(\sum_{i=1}^n\,\la_{ji}\log\rho_i)$ for $ j=1,\dots, n$.
It is obvious that $\vphi$ is a real analytic isomorphism
of $\BC^{*n}$ to $\BC^{*n}$.
Let us consider the Milnor fibrations of $f(\bfz,\bar\bfz)$ and 
$\hat f(\bfw)$ in the respective ambient tori $\BC^{*n}$. 
\[
 f:\BC^{*n}\backslash f\inv(0)\to \BC^*,\quad
\hat f:\BC^{*n}\backslash\hat f\inv(0)\to \BC^*
\]
Recall that the monodromy maps $h^*,\,\hat h^*$ are given
as 
\begin{eqnarray*}
&h^*:\,F^*\to F^*,\quad
\bfz\mapsto \exp(2\pi i/m_p)\circ \bfz\\
&\hat h^*:\,\hat F^*\to \hat F^*,\quad
\bfw\mapsto \exp(2\pi i/m_p)\circ \bfw.
\end{eqnarray*}
Recall that the $\BC^*$-action associated with $\hat f(\bfw)$ is the polar
 action of $f(\bfz,\bar\bfz)$.
Namely 
$\exp{i\theta}\circ \bfw=(\exp(ip_1\theta)w_1,\dots, \exp(ip_n\theta)w_n)$.
Thus we have the commutative diagram:
\begin{eqnarray*}
&\begin{matrix}
F_\al^*&\mapright{h^*}&F_\al^*\\
\mapdown{\vphi}&&\mapdown{\vphi}\\
\hat F_\al^*&\mapright{\hat h^*}&\hat F_\al^*
\end{matrix}\end{eqnarray*}
where 
$F_\al^*=f\inv(\al)\cap \BC^{*n}$ and 
$\hat F_\al^*={\hat f}\inv(\al)\cap \BC^{*n}$ for $\al\in \BC^*$.
\end{proof}
\subsubsection{Remark}
The case $f(\bfz,\bar\bfz)=z_1^{a_1}\bar
 z_1+\cdots+z_n^{a_n}\bar z_n$ is studied in \cite{R-S-V}.

\subsection{zeta-functions}
Now we know that  by \cite{OkaNagoya,Okabook},
the inclusion map
$\hat F^*\hookrightarrow \BC^{*n}$
is $(s-1)$-equivalence and
$\chi(\hat F^*)=(-1)^{n-1} \det (N-M)$ for $s=n$ and $0$ otherwise.
Note also the monodromy map
$\hat h:\hat F^*\to \hat F^*$ has a period $m_p$.
The fixed point locus of $(\hat h)^k$ is 
$F^*$ if  $m_p\,|\,k$ and $\emptyset$ otherwise.
Thus using the formula of the zeta function
(see, for example \cite{Milnor}),
\[
 \zeta_{\hat h^*}(t)=\exp(\sum_{j=0}^\infty\, (-1)^{n-1}d\, t^{jm_p}/(jm_p))=
(1-t^{m_p})^{(-1)^{n}d/m_p}
\]
where $d=\det\,(N-M)$
if $s=n$
and $d=0$ for $s<n$. Translating this in the monodromy $h^*:F^*\to F^*$, we
obtain
\begin{Corollary}\label{Lefschetz} $F^*$ has a homotopy type of CW-complex of dimension
 $n-1$
and the inclusion map
$F^*\hookrightarrow \BC^{*n}$ is an $(s-1)$-equivalence.
The zeta function $\zeta_{h^*}(t)$ of $h^*:F^*\to F^*$ is given
as $(1-t^{m_p})^{(-1)^{n}d/m_p}$ with  $d=\det\,(N-M)$
if $s=n$
and $\zeta_{h^*}(t)=1$ for $s<n$.
\end{Corollary}
\noindent
\subsubsection{ Remark}\label{restriction} In general, the restriction of the polar action
on $\BC^n$ to
$\BC^{*I}$ may not effective and to make the action effective,
we need to define polar weights as
$p_{I,i}=p_i/r_I$ and $m_{I,P}=m_p/r_I$ where $r_I$ is the gratest common
divisot of $\{p_i\,|\, i\in I\}$. Hoever the monodromy map
$h_I: F^{*I}\to F^{*I}$ is equal to the restriction of $h:F\to F$.

\section{Connectivity of $F$}
Now we are ready to patch together the information of the strata
$F^{*I}$ for the topology of $F$.
First we introduce the notion of {\em $k$-convenience}
which is introduced for holomorphic functions (\cite{Okabook}).
We say $f(\bfz,\bar\bfz)$ is {\em $k$-convenient} if
$f^I\not \cong 0$ for any $I\subset \{1,2,\dots, n\}$ with $|I|\ge  k$.
The following is obvious by the definition.
\begin{Proposition}
Assume that $f(\bfz,\bar\bfz)$ is a simple polar weighted homogeneous 
polynomial with $s$ monomials and assume that $f$ is
 $k$-convenient. Then 
$k\le s-1$.
\end{Proposition}
Now we have the following  result about the connectivity of $F$.
\begin{Theorem}\label{connectivity}
Assume that $f(\bfz,\bar\bfz)$ is a simple polar weighted homogeneous 
polynomial with $s$ monomials and assume that $f$ is
 $k$-convenient. Then $F$ is $\mini(k,n-2)$-connected.
\end{Theorem}
For the proof, we show  the following stronger assertion.
Let $I\subset \{1,2,\dots,n\}$ and
put
\begin{eqnarray*}
& \BC^n(*I)=\{\bfz=(z_1,\dots,z_n)\in \BC^n\,|\, z_j\ne 0,\,j\in I\},\,
F(*I)=F\cap \BC^n(*I).\\
&\BC^{*I}=\{\bfz\in \BC^n\,|\, z_j\ne 0\,  \text{iff}\,j\in I\},\,
F^{*I}=F\cap \BC^{*I}.
\end{eqnarray*}
\begin{Lemma}\label{stronger-connectivity}
Under the assumption as in Theorem \ref{connectivity},
the inclusion $F(*I)\hookrightarrow \BC^n(*I)$ is $\mini(k+1,n-1)$-equivalence.
\end{Lemma}
We prove the assertion on double induction on $(n,k)$.
Put 
\begin{eqnarray*}
& I_j=\{j,\dots, n\},\quad  K_j=\{1,\dots,\overset\vee j,\dots, n\}\\
&\BC^{n-1}_j=\BC^{K_j}=\BC^n\cap\{z_j=0\},\, F_j=F\cap\BC^{n-1}_j.
\end{eqnarray*}
Note that $F_j$ is 
the Milnor fiber of
$f^{K_j}$.
Theorem \ref{connectivity} follows from Lemma \ref{stronger-connectivity}
by taking 
$I=\emptyset$. Changing the ordering if necessary, we may assume that
$I=I_t$ for 
some $t$.
We consider the filtration of $F$:
\[
 F^*=F(*I_1)\subset F(*I_2)\subset F(*I_3) \subset \cdots\subset
 F(*I_{n})\subset F=F(*\emptyset).
\]
A key lemma is
\begin{Lemma}\label{k-eq}
The inclusion map $(F(*I_j),F(*I_{j-1}))\hookrightarrow (\BC^n(*I_j),
 \BC^n(*I_{j-1}))$ is $\mini(k+1,n-1)$-equivalence.
\end{Lemma}\begin{proof}

Let $T_j$ be a tubular neighborhood of $\{z_j=0\}$ in $\BC^n(*I_{j+1})$ 
such that
$T_j\cap F(*I_{j+1})$ is a tubular neighborhood of
	    $F_j(*I_{j+1})=\{z_j=0\}\cap  F(*I_{j+1})$ in
$ F(*I_{j+1}) $.
Consider the following diagram follows by the  excision
isomorphisms and  the Thom
	    isomorphisms $\psi$
for $D^2$-bundle:
\begin{multline*}
\begin{matrix}
H_{\ell+1} (F(*I_{j+1}),F(*I_{j}))&\mapright{\cong}&
H_{\ell+1} (F(*I_{j+1})\cap T_j,F(*I_{j})\cap T_j)\\
\mapdown{\tau_j}&&\mapdown{\tau_j'}&\\
H_{\ell+1} (\BC^n(*I_{j+1}),\BC^n(*I_{j}))&\mapright{\cong}&
H_{\ell+1} (T_j,\BC^n(*I_{j})\cap T_j)
\end{matrix}\\
\begin{matrix}
&\\
\mapright{\psi}&H_{\ell-1}(F_j(*I_{j+1}))\\
&\mapdown{\tau_j''}\\
\mapright{\psi}&H_{\ell-1}(\BC^{n-1}_j(*I_{j+1}))
\end{matrix}
\qquad
\end{multline*}
Now note that $f^{K_j}$ is ${(k-1)}$-convenient. Thus by the induction's
assumption on Lemma \ref{k-eq},  $\tau_j''$ is isomorphism for $\ell-1\le k-1 $.
This implies that 
$\tau_j',\, \tau_j$ is isomorphism for $\ell+1\le k+1$.
\end{proof}
{\em Proof of Lemma \ref{stronger-connectivity}.}
Now we can prove Lemma \ref{stronger-connectivity} by the induction on $j$ and Five Lemma,
assuming $I=I_{j}$ for some $j$,
applied to  two exact sequences for the pairs
$(F(*I_{j+1}),F(*I_{j}))$ and $(\BC^n(*I_{j+1}),F(*I_{j}))$ and
commutative diagrams:
\[
 \begin{matrix}
H_{\ell+1}(F(*I_{j+1}),F(*I_{j})))&\to&H_\ell(F(*I_{j}))&\to&
H_\ell(F(*I_{j+1}))\\
\mapdown{\tau_j}&&\mapdown{\iota_{j}}&&\mapdown{\iota_j}\\
H_{\ell+1}(\BC^n(*I_{j+1}),\BC^n(*I_{j})))&\to&H_\ell(\BC^n(*I_{j}))&\to&
H_\ell(\BC^n(*I_{j+1}))
\end{matrix}
\]
Induction starts for $j=1$: $\iota_1$ is $min(k+1,n-1)$-equivalence
by Corollary \ref{Lefschetz}.
This completes the proof of Lemma \ref{stronger-connectivity}.
\qed
\subsection{Euler numbers and zeta functions}
Let $f(\bfz,\bar\bfz)=\sum_{j=1}^s c_j\bfz^{\bfn_j}\bar\bfz^{\bfm_j}$
 be a simplicial polar weighted homogeneous. Let
$$\cS=\{I\subset \{1,\dots,n\}; f^I \,\text{is}\, \,full\}$$
and put $r_I=\gcd_{i\in I}\{p_i\}$
and $m_{p,I}=m_p/r_I$ and put $d_I=|\det_{i\in I}(\bfn_i-\bfm_i)|$.
Thus for $I\in \cS$, $f^I$ is a simplicial full polar weighted
 homogeneous polynomial of polar weight type
$(p_i/r_I)_{i\in I}$ with degree $m_{p,I}$. 
We observed in Remark \ref{restriction} that
the monodromy map $h^{*I}:F^{*I}\to F^{*I}$ is equal to the restriction
of the monodromy map $h:F\to F$.
We denote the zeta function
of the monodromy map
\[
 h:F\to F,\, h^{*I}=h|_{F^{*I}}:F^{*I}\to F^{*I}
\]
by $\zeta(t),\,\zeta^{*I}(t)$ respectively.
Recall  that $\zeta(t)$ is a  alternating product of characteristic
polynomials(\cite{Milnor}). Namely  
\[
 \zeta(t)=\prod_{j=0}^{n-1}P_j(t)^{(-1)^{j+1}}
\]
where $P_j$ is the characteristic polynomial of the monodromy action
on $h_*:H_j(F,\BQ)\to H_j(F,\BQ) $.
By Theorem \ref{main-result} and
the additive formula for the Euler characteristics, using a similar
argument as that of Proposition 2.8, \cite{Okabook}, we have:
\begin{Theorem}
\begin{enumerate}
\item $\chi(F)=\sum_{I\in \cS}(-1)^{|I|-1}d_I$.
\item
$\zeta(t)=\prod_{I\in \cS}\zeta^{*I}(t),\,
\zeta^{*I}(t)=(1-t^{m_{p,I}})^{(-1)^{|I|}d_I/m_{p,I}}$.
\end{enumerate}
\end{Theorem}
\subsection{Examples}

1. Assume that 
$f_1(\bfz)$ is a homogeneous polynomial defined by 
\[
 f_1(\bfz)=z_1^{a_1}+z_2^{a_2}+\cdots+z_n^{a_n},\quad a_1,\dots,a_n\ge 2.
\]
Then $F=f_1\inv(1)$ is $(n-2)$-connected and 
\[
 \chi(F)=\sum_{j=1}^n
\sum_{|I|=j}\chi(F^{*I})=(a_1-1)(a_2-1)\cdots(a_n-1)-(-1)^n\]
and 
\[
 {\rm div}( \zeta_h)=(\La_{a_1}-1)\cdots (\La_{a_n}-1)-(-1)^n
\]
as is well-known by \cite{Pham,Brieskorn:1966,MilnorOrlik}.
Here ${\rm div}((t-\la_1)\cdots(t-\la_k))=\sum_{i=1}^k \la_i\in \BZ\cdot
\BC^*$
and $\La_m={\rm div}(t^m-1)$.

2. Consider
\[
 f_2(\bfz,\bar\bfz)=z_1^{a_1}\bar z_2+\cdots+z_{n-1}^{a_{n-1}}\bar z_n+z_n^{a_n}
\]
Then $f_2$ is a simplicial polar weighted polynomial.
and 
\[
 \cS=\{I_j=\{1,\dots,j\}\,|\,j=0,\dots,n-1\}.
\]
Thus we have
\begin{eqnarray*}
&\chi(F)=(-1)^{n-1}\left(
a_1a_2\cdots a_n-a_2\cdots a_n+\cdots+(-1)^{n-1}a_n
\right)\\
&\log\zeta(t)=(-1)^{n}\left(
\frac{1}{(1-t^{a_1\cdots a_n})}-\frac 1{(1-t^{a_2\cdots a_n})}+
\cdots+(-1)^{n-1}\frac 1{(1-t^{a_n})}
\right)
\end{eqnarray*}
\begin{proof}
The polar weight of $f_2$ is given by $(p_1,\dots, p_n;m_p)$
where 
\begin{eqnarray*}
& m_p=a_1\cdots a_n,\, p_1=m_p\left(\frac 1{a_1}+\cdots +\frac{1}{a_1\cdots
 a_n}\right),\\
&p_2=m_p\left(\frac 1{a_2}+\cdots+\frac 1{a_2\cdots a_n}\right)\\
&\qquad \vdots\\
&p_{n-1}=m_p\left(\frac 1{a_{n-1}}+\frac 1{a_{n-1}a_n}\right),\,
 p_n=\frac {m_p}{a_n}
\end{eqnarray*}
Thus the assertion follows from Corollary \ref{Lefschetz}.
\end{proof}
\subsection{Surface cases} Consider the case $n=3$. We consider two
simplicial polar weighted 
homogeneous  polynomials.
\begin{eqnarray*}
&f_1(\bfz,\bar\bfz)=z_1^{a_1}\bar z_2^{b_1}+z_2^{a_2}\bar
 z_3^{b_2}+z_3^{a_3},\,\,a_1,a_2, b_1,b_2>0\\
&f_2(\bfz,\bar\bfz)=z_1^{a_1}\bar z_2^{b_1}+z_2^{a_2}\bar
 z_3^{b_2}+z_3^{a_3}\bar z_1^{b_3},\,\,a_1a_2a_3>b_1b_2b_3>0.
\end{eqnarray*}
They are 1-convenient.
Let $F_1=f_1\inv(1)$ and $F_2=f_2\inv(1)$. By Theorem \ref{connectivity},
$F_1,\,F_2$ are simply connected.
Their Betti numbers $b_2(F_i)$ are given as
\[
 b_2(F_1)=a_1a_2a_3-a_2a_3+a_3-1,\quad
b_2(F_2)=a_1a_2a_3-b_1b_2b_3-1.
\]
(I) First we consider $f_1$. The normalized polar weight for $f_1$ is given as
\begin{eqnarray*}
v_1=\frac{b_1b_2}{a_1a_2a_3}+\frac{b_1}{a_1a_2}+\frac 1{a_1},
\,v_2=\frac {b_2}{a_2a_3}+\frac 1{a_2},\, v_3=\frac{1}{a_3}
\end{eqnarray*} 
 Let $ r=\gcd(b_1b_2,a_1a_2a_3),\, r_1=\gcd(b_2,a_2a_3)$.
Then $m_p$ is given as $a_1a_2a_3/r$
and the zeta function of $h_1:F_1\to F_1$ is given as
\begin{eqnarray*}
&\zeta_{h_1}(t)=P_{0}(t)^{-1}P_{2}(t)^{-1}=\frac
{(1-t^{a_2a_3/r_1})^{r_1}}{(1-t^{a_1a_2a_3/r})^r(1-t^{a_3})}\\
\end{eqnarray*}
where $P_{2}(t)$ is the characteristic polynomial of the monodromy action
 $h_{1*}:H_2(F_1;\BQ)\to H_2(F_1;\BQ)$. Note that $P_{0}(t)=1-t$.
For example, 
\begin{eqnarray*}
&\zeta_{h_1}(t)=
\frac
{(1-t^{a_2a_3})}{(1-t^{a_1a_2a_3})(1-t^{a_3})},\,\, b_1=b_2=1\\
&\zeta_{h_1}(t)=\frac
{(1-t^{a_2'a_3})^{2}}{(1-t^{a_1'a_2'a_3})^{4}(1-t^{a_3})},\,\, a_1=2a_1',\,a_2=2a_2',\,b_1=b_2=2.\\
\end{eqnarray*}
(II) We consider $f_2$.  The normalized polar weight for $f_2$ is
given as:
\begin{eqnarray*}
&v_1=\frac{a_2a_3+b_1a_3+b_1b_2}{a_1a_2a_3-b_1b_2b_3},
\, v_2=\frac{a_1a_3+a_1b_2+b_2b_3}{a_1a_2a_3-b_1b_2b_3},
\,\,v_3=\frac{a_1a_2+a_2b_3+b_1b_3}{a_1a_2a_3-b_1b_2b_3}.
\end{eqnarray*}
Put $d=a_1a_2a_3-b_1b_2b_3$.
The least common multiple $m_p$ of the denominators of $v_1,v_2,v_3$
depends on $\gcd(d,a_2a_3+b_1a_3+b_1b_2)$
and so on. We only gives two examples.

(1) Assume that $a_1=a_2=a_3=a,\,b_1=b_2=b_3=b$. Then
$v_1=v_2=v_3=\frac 1{a-b}$. Thus
\[
 \zeta_{h_2}(t)=(1-t^{a-b})^{a^2+ab+b^2}.
\]

(2) Assume that
$\gcd(d,a_2a_3+b_1a_3+b_1b_2)=\gcd(d,a_1a_3+a_1b_2+b_2b_3)=$
$\gcd(d,a_1a_2+a_2b_3+b_1b_3)=1$. Then $m_p=d$  and
$ \zeta_{h_2}(t)=(1-t^d)$.

For example, if $a_1=2,a_2=3,a_3=5$ and $b_1=b_2=b_3=1$,
we get
$\zeta_{h_2}(t)=(1-t^{29})$.

\begin{thebibliography}{10}

\bibitem{Brieskorn:1966}
E.~Brieskorn.
\newblock Beispiele zur {D}ifferentialtopologie von {S}ingularit\"aten.
\newblock {\em Invent. Math.}, 2:1--14, 1966.

\bibitem{Molina}
J.~Cisneros-Molina.
\newblock Join theorem for polar weighted homogeneous singularities.
\newblock In {\em Proceeding of Le-Fest, Cuernavaca, 2007, to appear}.

\bibitem{H-Lopez}
L.~Hern{\'a}ndez de~la Cruz and S.~L{\'o}pez~de Medrano.
\newblock Some families of isolated singularities.
\newblock In {\em Proceeding of Le-Fest, Cuernavaca, 2007, to appear}.

\bibitem{Milnor}
J.~Milnor.
\newblock {\em Singular Points of Complex Hypersurface}, volume~61 of {\em
  Annals Math. Studies}.
\newblock Princeton Univ. Press, 1968.

\bibitem{MilnorOrlik}
J.~Milnor and P.~Orlik.
\newblock Isolated singularities defined by weighted homogeneous polynomials.
\newblock {\em Topology}, 9:385--393, 1970.

\bibitem{OkaJoin}
M.~Oka.
\newblock On the homotopy types of hypersurfaces defined by weighted
  homogeneous polynomials.
\newblock {\em Topology}, 12:19--32, 1973.

\bibitem{OkaNagoya}
M.~Oka.
\newblock On the topology of full nondegenerate complete intersection variety.
\newblock {\em Nagoya Math. J.}, 121:137--148, 1991.

\bibitem{Okabook}
M.~Oka.
\newblock {\em Non-degenerate complete intersection singularity}.
\newblock Hermann, Paris, 1997.

\bibitem{Pham}
F.~Pham.
\newblock Formules de {P}icard-{L}efschetz g\'en\'eralis\'ees et ramification
  des int\'egrales.
\newblock {\em Bull. Soc. Math. France}, 93:333--367, 1965.

\bibitem{Pichon-Seade}
A.~Pichon and J.~Seade.
\newblock Real singularities and open-book decompositions of the 3-sphere.
\newblock {\em Ann. Fac. Sci. Toulouse Math. (6)}, 12(2):245--265, 2003.

\bibitem{Pichon-Seade2}
A.~Pichon and J.~Seade.
\newblock Fibered multilinks and singularities $f{\bar g}^*$.
\newblock {\em Preprint, Nov. 13}, 2007.

\bibitem{R-S-V}
M.~A.~S. Ruas, J.~Seade, and A.~Verjovsky.
\newblock On real singularities with a {M}ilnor fibration.
\newblock In {\em Trends in singularities}, Trends Math., pages 191--213.
  Birkh\"auser, Basel, 2002.

\bibitem{Seade}
J.~Seade.
\newblock Open book decompositions associated to holomorphic vector fields.
\newblock {\em Bol. Soc. Mat. Mexicana (3)}, 3(2):323--335, 1997.

\bibitem{SeadeBook}
J.~Seade.
\newblock On the topology of hypersurface singularities.
\newblock In {\em Real and complex singularities}, volume 232 of {\em Lecture
  Notes in Pure and Appl. Math.}, pages 201--205. Dekker, New York, 2003.

\end{thebibliography}
\def\cprime{$'$} \def\cprime{$'$} \def\cprime{$'$} \def\cprime{$'$}
  \def\cprime{$'$}

\end{document}